\documentclass[journal,onecolumn]{IEEEtran}
\usepackage{amsmath,epsfig}
\usepackage{graphicx}
\usepackage{multirow}
\usepackage{threeparttable}
\usepackage{amsfonts}
\usepackage{lscape}
\usepackage{color}
\usepackage{threeparttable}
\usepackage{tabularx}
\usepackage{amsfonts}
\usepackage{algpseudocode}
\usepackage{caption}
\usepackage{algorithm}
\usepackage{amsmath,epsfig}
\usepackage{graphicx}
\usepackage{multirow}
\usepackage{threeparttable}
\usepackage{amsfonts}
\usepackage{lscape}
\usepackage{amsmath,graphicx}
\usepackage{color}
\usepackage{threeparttable}
\usepackage{tabularx}
\usepackage{amsfonts}
\usepackage{algpseudocode}
\usepackage{caption}
\usepackage{algorithm}
\usepackage{subcaption}
\usepackage{makecell}
\usepackage{hyperref}
\hypersetup{
    colorlinks,
    citecolor=black,
    filecolor=black,
    linkcolor=black,
    urlcolor=black
}

\ifCLASSINFOpdf
\else
\fi
\begin{document}
%
\title{On Algorithms for Solving the Rubik's Cube}
%
%
%

\author{Ahmad~Kaleem
        and Ahsan~Kaleem
}

\maketitle

\begin{abstract}
In this paper, we present a novel algorithm and its three variations for solving the Rubik's cube more efficiently. This algorithm can be used to solve the complete $n \times n \times n$ cube in $O(\frac{n^2}{\log n})$ moves. This algorithm can also be useful in certain cases for speedcubers. We will prove that our algorithm always works and then perform a basic analysis on the algorithm to determine its algorithmic complexity of $O(n^2)$. Finally, we further optimize this complexity to $O(\frac{n^2}{\log n})$.
\end{abstract}


%
\IEEEpeerreviewmaketitle

\section{Introduction}
\label{sec:intro}

An $n \times n \times n$ Rubik's cube is composed of $n^3$ \textbf{cubies} ($1 \times 1 \times 1$ cubes) each of which is located at a position $(i,j,k)$ where $i,j,k \in \{0,1,...,n-1\}$. Each cubie in a Rubik's cube has a color on a face where there are a total of six faces on the cube with six unique colors. The six faces can be denoted by up, down, right, left, front and back which represents the position of the face with respect to the viewing angle. If the cube is placed right in front of you, the face right in front of you will be the front face and so on. A Rubik's cube is solved when every cubie on each face is the same color as the rest of the cubies on that same face (i.e. each face is the same color). A \textbf{slice} on a Rubik's cube is a set of cubies with a common $x$, $y$ or $z$ coordinate. A \textbf{legal move} consists of rotating one slice of the Rubik's cube 90 degrees in either direction. A \textbf{cubie cluster} is a set of cubies on the cube such that any sequence of legal moves preserves the contents of the cluster (i.e. any cubie originally in the cluster always remains in the cluster while any cubie not in some cluster can never come to that specific cluster). There are 3 types of cubies on a Rubik's cube: center cubies, edge cubies and corner cubies. A corner cubie is a cubie on the corner of a cube which has 3 different visible faces with three colors. An edge cubie is a cubie on the intersection of two faces of the cube and thus has two visible colors. Finally a center cubie is one in the center of the cube which is on only one face and thus has only one visible face. To refer to a specific corner cubie we use three letters which refer to the three faces on which this piece is located. As an example $UFR$ refers to the corner cubie located at the intersection of the up, front and right face. To refer to a edge cubie, we use the two faces which it is on and the slice which it is on. For example, $UFV_i$ refers to the cubie on the up and front and the $i$th vertical slice. This also corresponds to the notation for moves which are mentioned below. To refer to a single center cubie, we use one face letter and two slices. For example $FV_iH_j$ refers to the cubie on the front face and the intersection of the $i$th vertical and $j$th horizontal slice. We also define the notation we will be using. Regular Rubik's cube notation for face moves with letters representing each face is used along with more general notation in terms of slices is used. Both notations are shown more clearly in Fig. \ref{fig:fig1} and Fig. \ref{fig:fig2} below.

\begin{figure}[hbt!]
    \includegraphics[width=15cm]{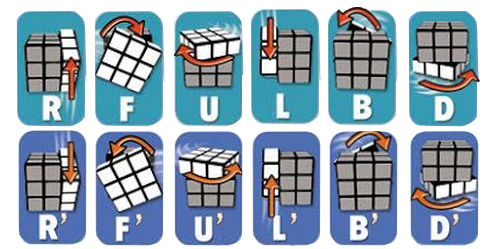}
    \caption{Notation for face moves}
    \label{fig:fig1}
\end{figure}

\begin{figure}[hbt!]
    \includegraphics[width=15cm]{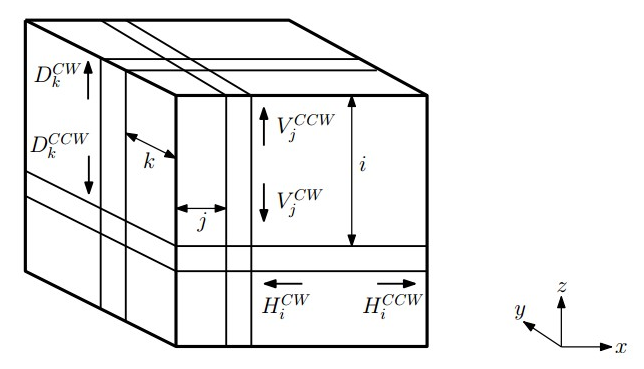}
    \caption{General notation for slice moves \cite{Ref1}}
    \label{fig:fig2}
\end{figure}

\section{Algorithm}\label{sec:alg}

\subsection{Case 1 (Face Moves)}\label{ssec:case1}

\begin{figure}[hbt!]
    \includegraphics[width=7cm]{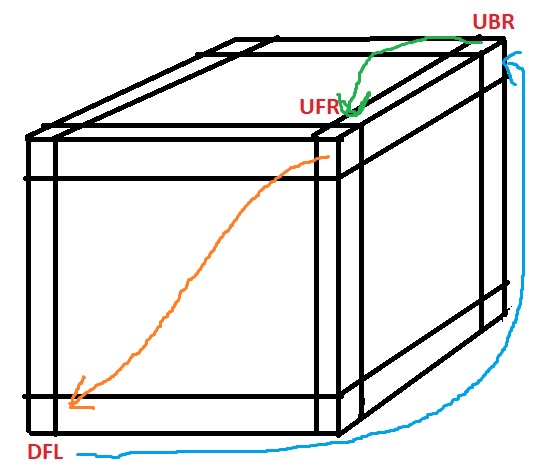}
    \caption{Permutation Cycle in Case 1}
    \label{fig:fig3}
\end{figure}

In the first case of our algorithm we use only face moves. We claim that the algorithm $U\ R'\ D\ R\ U'\ R'\ D'\ R$ swaps three corner cubies while applying the identity permutation to all other cubie. More specifically, the corners $DFL$, $UBR$ and $UFR$ are swapped in that order i.e. $DFL$ goes to $UBR$, $UBR$ goes to $UFR$ and $UFR$ goes to $DFL$ while all other cubies remain unaffected. It can also be represented by the permutation cycle: ($DFL$ $UBR$ $UFR$). This can be seen more clearly in Fig. \ref{fig:fig3}. 

\textit{Proof.}
We clearly do not need to consider the cubies which are not affected by any of the $U$, $R$ or $D$ type moves since they will remain unaffected throughout the algorithm. We thus start by considering the three corner cubies to show that they are indeed cycled as desired. First consider the corner $DFL$. The first move which affects it is $D$ which moves it to $DFR$. Then $R$ takes it to $UFR$, $U'$ takes it to $UBR$, $R'$ brings it back to $UFR$ and then finally $R$ takes it to $UBR$. Next, consider the corner $UBR$. It first goes to $UFR$ with a $U$ move, then to $DFR$ with $R'$, $DBR$ with $D$ then back to $DFR$ with a $R$, back to $DBR$ with $R'$, to $DFR$ with a $D'$ and then finally to $UFR$ with an $R$. Finally consider the corner $UFR$ as the algorithm is applied. It first goes to $UFL$ with a $U$ move, then back to $UFR$ with $U'$, down to $DFR$ with an $R'$ and then finally to $DFL$ with a $D'$. We now prove that no other cubies are affected. Since in the algorithm, each move has its inverse present exactly once (e.g. for the $U$ a $U'$ move exists), the net effect on the remaining cubies is that of the identity move. It can be manually verified that each of the cubies other than $DFL$, $UBR$ and $UFR$ which are affected by at least one of the moves will be brought back to their original locations by the algorithm since each move will be canceled by its inverse. Since the size of the cube does not affect the algorithm (as only face moves are used) this algorithm will hold true for any $n \times n \times n$ cube by testing on any one size such as a $3 \times 3 \times 3$ which we have verified. $\Box$

\subsection{Case 2 (Inner Layers)}\label{ssec:case2}
In the second case, we extend our algorithm to use inner layers. More specifically, we claim that the general algorithm $H_0^{CW}$ o  $V_x^{CW}$ o $H_{n-1}^{CCW}$ o $V_x^{CCW}$ o $H_0^{CCW}$ o $V_x^{CW}$ o $H_{n-1}^{CW}$ o $V_x^{CCW}$  permutes three cubies while applying the identity permutation to the rest of the cubies. Here $x$ can be any number from 0 to $n-1$. In other words, any vertical layer can be rotated throughout the algorithm as long as it remains the same throughout. This time the cubies which are cycled are $DLD_{n-x-1}$, $URD_x$ and $UFV_x$. The permutation cycle is ($DLD_{n-x-1}$ $URD_x$ $UFV_x$). So if $x= n-2$, the cubies which are permuted are edge cubies to the left by one piece relative to $DFL$, $UBR$ and $UFR$. For simplicity in the proof below, we refer to these cubies as $DFL'$, $UBR'$ and $UFR'$ respectively.

\textit{Proof.}
We will only consider the cases when $x$ is not equal to $0$ or $n-1$ since in this case all the moves are face moves which we have already considered in the Case 1. We first prove that the cubies which are permuted are indeed permuted in such a way. First consider cubie $DFL'$. It is unaffected by the first two moves, then it goes to the front face, then the top, then to the top-right which is $UBR'$ in this case and then remains there. Next consider the cubie $UBR'$. It first comes to the front face, then to the bottom, then to the right, then remains unaffected for the next three moves, then comes to the front and then to $UFR'$ where it ends up. Finally consider the cubie $UFR'$. It first goes to the left face, then remains unaffected for the next three moves, then back to the front face, then to the bottom and then to $DFL'$ where it ends up. Next we prove that all other cubies are unaffected by this algorithm. Clearly the ones which are not moved at all by any of those moves are unaffected.

We thus consider the cubies which are affected by at least one of the moves. First consider a cubie which is affected by only one type of move (e.g. the center pieces). Since the algorithm has an inverse move present for every move, the net effect of applying the algorithm will be the identity permutation since only one type of move ever affects it and that type of move gets reversed by its inverse which is present as we showed in case 1. For example if a cubie is only affected by $V_x$ moves, the algorithm becomes $V_x^{CW}$ o $V_x^{CCW}$ o $V_x^{CW}$ o $V_x^{CCW}$, $V_x^{CW}$ o $V_x^{CCW}$ or $V_x^{CCW}$ o $V_x^{CW}$ which are all equal to the identity permutation since the moves are canceled out by their inverse.

Next consider the cubies affected by exactly two types of moves. We consider all possible cases. First let a cubie be affected by only $H_0$ and $V_x$ moves. We consider all possible subsets of the $H_0$ and $V_x$ moves which such a cubie can be affected by. First let it be affected by all of the $H_0$ and $V_x$ moves in the algorithm. Then clearly it must originally be in the top layer since the first move in the algorithm is $H_0$. Moreover it must be in the layer originally on the top face which aligns with the $V_x$ layer after the $H_0$ move at the beginning of the algorithm. The next move ($V_x^{CW}$) brings this layer to the front face and then the ones on this layer affected by two types of moves go back to the top layer with $V_x^{CCW}$. After that the cubies in this layer are rotated counterclockwise with $H_0^{CCW}$  and then there is only one cubie from that original layer left which aligns with the $V_x$ layer (1) so this is the only cubie which will be affected by the next $V_x$ moves. Thus there is only one cubie which is affected by all of the $H_0$ and $V_x$ moves in the algorithm. For this cubie, the algorithm simplifies to $H_0^{CW}$ o $V_x^{CW}$ o $V_x^{CCW}$ o $H_0^{CCW}$ o $V_x^{CW}$ o $V_x^{CCW}$ which has no effect since consecutive terms cancel out and the result is the identity permutation. Now if the cubie was not aligned with the $V_x$ layer at (1), then the algorithm applied to it will be $H_0^{CW}$ o $V_x^{CW}$ o $V_x^{CCW}$ o $H_0^{CCW}$ which again simplifies to the identity permutation. Now consider cubies which are not affected by the first $H_0$ move. Clearly such a cubie is not in the top layer originally. Since we are considering the cubies which are affected by both $V_x$ and $H_0$ moves, to be affected by $H_0$ moves at any time, the cubie must come in the top face at some point. However the only cubies that can come are from the back face (with the $V_x^{CW}$ move after the first $H_0$ move) but since a $V_x^{CCW}$ move comes before the next $H_0$ move, these cubies are not affected by $H_0$ moves and so are part of the first case in which a cubie is affected by only one type of move. The reason that no other cubie can come to the top face is because a $V_x^{CCW}$ precedes a $V_x^{CW}$ move at both points. Thus it is not possible for a cubie to be affected by both $V_x$ and $H_0$ moves without being affected by the first $H_0$ move in the algorithm.

	Now we will look at the case of cubies being affected by $V_x$ and $H_{n-1}$ moves only. Using a similar analysis as for $V_x$ and $H_0$ we find that the subset of the algorithms for these cubies can be one of the following two: $V_x^{CW}$ o $H_{n-1}^{CCW}$ o $V_x^{CCW}$ o $V_x^{CW}$ o $H_{n-1}^{CW}$ o $V_x^{CCW}$ or $V_x^{CW}$ o $H_{n-1}^{CCW}$ o $H_{n-1}^{CW}$ o $V_x^{CCW}$. Again, in both of these cases the moves get canceled by their reverse counterpart, and hence become equivalent to the identity permutation.

 Finally let a cubie be affected by only $H_0$ and $H_{n-1}$ moves. Since these moves are disjoint and cannot overlap, it is not possible to have any cubie be affected by only these moves unless it is affected by a $V_x$ move to bring it on the bottom face in the first place which is a contradiction since we assumed that such a cubie is affected by only two of the three possible moves.

We now show that the only cubies which are affected by all three types of moves are the three cubies which are permuted and we have already shown that the algorithm correctly permutes them. This is because if we consider the cubies on the top face for example, only one cubie at a time can be brought to the bottom face since we do not have two consecutive $V_x$ moves in the same direction at any point in the algorithm. Moreover since we only have one $H_0$ move in each direction, the maximum number of cubies which can be brought to the bottom face is two and these are indeed the ones which our algorithm affects. Clearly if a cubie is not brought to the bottom face, it cannot be affected by all three types of moves since the $H_{n-1}$ move only affects the bottom face. Thus the only cubies which are affected by all three types of moves are the ones which are permuted. Thus these cubies are permuted while all remaining cubies are left unaffected as desired. $\Box$

\subsection{Case 3 (Middle Slices)}\label{ssec:case3}
In the third case of our algorithm, we use middle slices to permute three center cubies. This algorithm exists for $n \times n \times n$ cubes which have $n > 3$. It doesn't work for $n \leq 3$ since there is only one center slice and so there is only one move possible for the middle slice. The algorithm is as follows: $H_k^{CW}$ o $V_{n-1}^{CW}$ o $H_{n-k-1}^{CCW}$ o $V_{n-1}^{CCW}$ o $H_k^{CCW}$ o $V_{n-1}^{CW}$ o $H_{n-k-1}^{CW}$ o $V_{n-1}^{CCW}$. Here $k$ must be less than or equal to $\lfloor \frac{n}{2} \rfloor$ so that the $H_k$ move is on the upper half of the cube. The three cubies permuted are $FV_kH_{n-k-1}$, $RH_kV_k$ and $BH_kV_k$. As a permutation cycle, this is ($FV_kH_{n-k-1}$  $BH_kV_k$  $RH_kV_k$).

\textit{Proof.}
We first prove that these cubies are permuted as stated. First consider the cubie originally on the front face ($FV_kH_{n-k-1}$). It goes to the right face with the third move in the algorithm, then to the $H_k, V_k$ position on the right face with the $V_{n-1}$ move, then to the back face with the $H_k$ move where it stays and this is indeed the position it was supposed to go to. Next consider the cubie originally on the right face ($RH_kV_k$). With the first $H_k$ move, it comes to the front face, the next three moves do not affect it, with the $H_k^{CCW}$ move it goes to the right face, the $V_{n-1}^{CW}$ move rotates it so that it is in the $H_{n-k-1}$, $V_k$ position, then with the $H_{n-k-1}^{CW}$ move it returns to the front face in the correct position where it belongs. Finally consider the cubie originally on the back face ($BH_kV_k$). With the first $H_k$ move, it goes to the right face, then with the $V_{n-1}$ move, it is shifted to the $H_{n-k-1}, V_k$ position, with a $H_{n-k-1}^{CCW}$ move it goes to the back face, the next three moves do not affect it, then the $H_{n-k-1}^{CW}$ move brings it to the right face and finally with the $V_{n-1}^{CCW}$ move it goes to the $H_k, V_k$ position where it belongs.

	We now prove that all other cubies are unaffected by this algorithm. Again we ignore the cubies affected by none of the moves since those are clearly unaffected. We start by considering the cubies affected by only one of the three possible moves. In this case, the possible algorithms applied to a cubie are $H_k^{CW}$ o $H_k^{CCW}$, $H_{n-k-1}^{CCW}$ o $H_{n-k-1}^{CW}$, $V_{n-1}^{CW}$ o $V_{n-1}^{CCW}$, $V_{n-1}^{CCW}$ o $V_{n-1}^{CW}$ or $V_{n-1}^{CW}$ o $V_{n-1}^{CCW}$ o $V_{n-1}^{CW}$ o $V_{n-1}^{CCW}$. It is clear that these are the only possibilities since if a cubie is not affected by the inverse move, that means that it was moved away by another move to a different layer but that is a contradiction since we assumed that the cubies are affected by only one type of move. Furthermore, it can be seen that all the algorithms above are equivalent to the identity permutation.

	We now consider the cubies affected by $H_k$ and $V_{n-1}$ moves only. First let a cubie be affected by all of the $H_k$ and $V_{n-1}$ moves in the algorithm. Then the algorithm applied to this cubie simplifies to $H_k^{CW}$ o $V_{n-1}^{CW}$ o $V_{n-1}^{CCW}$ o $H_k^{CCW}$ o $V_{n-1}^{CW}$ o $V_{n-1}^{CCW}$  which is the identity permutation since each move is canceled out by its inverse. We now consider all possible subsets of this algorithm which can be applied to a cubie affected by only $H_k$ and $V_{n-1}$ moves. Since there are only two $H_k$ moves in the whole algorithm and the cubie must be affected by at least one of them, we first consider the cubies affected by the first $H_k$ move. Because of this, the cubie must originally be anywhere in the $H_k$ slice. After this move, we first let the cubie be affected by the next move (1). Then clearly it must also be affected by the next $V_{n-1}$ move since the cubies affected by the first two moves are currently in this slice. Next, such a cubie will be affected by the $H_k^{CCW}$ since it was affected the $H_k^{CW}$ move at the beginning of the algorithm and so it is at the same position right now. It will also be affected by the next $V_{n-1}^{CCW}$ move since at this point it is on the right face and then also the last move in the algorithm, $V_{n-1}^{CCW}$. But this is the same case as we considered before since the cubie was affected by all of the $V_{n-1}$ and $H_k$ moves. The other case at (1) is if it was not affected by the $V_{n-1}$ move. Then after this it will not affected by the next $V_{n-1}$ move in the algorithm since it is not in that slice but will be affected by the next $H_k^{CCW}$ and then both $V_{n-1}$ moves since it is in that slice. Thus if a cubie was affected by the first $H_k$ move, it will either be affected by all $H_k$ and $V_{n-1}$ moves or the algorithm $H_k^{CW}$ o $H_k^{CCW}$ o $V_{n-1}^{CW}$ o $V_{n-1}^{CCW}$ which also cancels out to give the identity permutation. Next we consider cubies affected by $V_{n-1}$ and $H_k$ type moves which are not affected by the first $H_k$ move at the beginning of the algorithm. Then clearly such a cubie must be affected by the second $H_k$ move in the algorithm since there are only two $H_k$ moves in the whole algorithm. But since it was not affected by the first $H_k$ move, it cannot be in the $H_k$ slice at the beginning of the algorithm. But this implies that such a cubie is affected only by $V_{n-1}$ moves since  the $V_{n-1}$ moves only rotate the right slice by a maximum of 90 degrees and so a cubie in a different horizontal layer at the beginning cannot coincide with the $H_k$ layer later on. But this is a contradiction since we are only considering the cubies affected by both $H_k$ and $V_{n-1}$ moves.

	Using a similar analysis, two possible algorithms can be found for the cubies affected by $V_{n-1}$ and $H_{n-k-1}$ moves only. One is the algorithm  $V_{n-1}^{CW}$ o $H_{n-k-1}^{CCW}$ o $V_{n-1}^{CCW}$ o $V_{n-1}^{CW}$ o $H_{n-k-1}^{CW}$ o $V_{n-1}^{CCW}$ in which case a cubie is affected by all of the $V_{n-1}$ and $H_{n-k-1}$ moves in the algorithm. Clearly this algorithm simplifies to the identity permutation since moves are canceled out by their inverse. In the other case where not all of the $V_{n-1}$ and $H_{n-k-1}$ moves are used, the algorithm applied becomes  $H_{n-k-1}^{CCW}$ o $V_{n-1}^{CCW}$ o $V_{n-1}^{CW}$ o $H_{n-k-1}^{CW}$ which again simplifies to the identity permutation.

	Finally, a cubie cannot be affected by $H_k$ and $H_{n-k-1}$ moves only since these two layers are disjoint and the only way a cubie can go between these layers is if a $V_{n-1}$ move is used but this is not possible in this case.

We now show that the only cubies which are affected by all three moves are the three cubies which are permuted and we have already shown that the algorithm correctly permutes them. This is because if we consider the cubies on the top face for example, only one cubie at a time can be brought to the bottom slice since we do not have two consecutive $V_{n-1}$ moves in the same direction at any point in the algorithm. Moreover since we only have one $H_0$ move in each direction, the maximum number of cubies which can be brought to the bottom slice is two and these are indeed the ones which our algorithm affects. Clearly if a cubie is not brought to the bottom slice, it cannot be affected by all three types of moves since the $H_{n-k-1}$ move only affects the bottom slice. Thus the only cubies which are affected by all three types of moves are the ones which are permuted. Thus these cubies are the only ones which are permuted while all remaining cubies are unaffected as desired. $\Box$

\subsection{Additional Remarks}\label{ssec:remarks}
We now have three algorithms which can permute any three cubies (corners, edges and centers) while applying the identity permutation to all other cubies in the cube. Using combinations of these algorithms for specific cubies, we can solve any cluster configuration with even parity in a cluster move solution of length $O(1)$. This is because by our previous proofs, we know that there exist a set of permutations which can be applied to a single cubie cluster while applying the identity permutation to every other cluster. Thus any even permutation can be written as a composition of these permutations and has an inverse that can also be written as the composition of these permutations. However, since each cluster has finite size, the inverse composition must have $O(1)$ length. Thus there exists a $O(1)$ length sequence of moves that can be applied to the cube that results in one cluster being solved with all other cubie clusters remaining unaffected. We call this \textbf{Lemma 1}.

\textbf{Note:} Our algorithm permutes cubies in the same cubie cluster. However, in the current version, if we consider four slices for the permuted cubies (up, down, right and left), the algorithm can only permute cubies between the up and down slices or the right and left slices for all three cases mentioned above. However if we want to permute cubies between, for example the up and right slices, the algorithm can be modified slightly by adding a move in the beginning and then adding the inverse of that move at the end of the algorithm. For example we can do a $D_n^{CW}$ move at the beginning of the algorithm in the second case and add a $D_n^{CCW}$ move at the end and this will permute cubies between the top slice and the right slice. Again all other cubies will remain unaffected since the moves added at the beginning will get canceled at the end by the inverse.

\section{Solving a Rubik's Cube naively using this algorithm}\label{sec:naiveSol}
In this paper, we have shown three variations of an algorithm which can together be used to solve all types of cubie clusters (corner, center and edge) by in turn permuting 3 cubies in the same cluster. If we use this without any optimization, we can always solve a $n\times n \times n$ cube in $O(n^2)$ moves.

To prove this we start by solving the corner cluster. Since there are always 8 corner cubies for any sized Rubik's cube, this step can be done in constant time with respect to n ($O(1)$ moves).

Moving on, we solve the edge cubie clusters. There are $4\cdot(n-2)$ edge cubies on one face of an $n\times n \times n$ cube. The number of unique edge clusters is $\left \lceil{\frac{n}{2}}\right \rceil -1$ since there are $\left \lceil{\frac{n}{2}}\right \rceil$ pieces on one side of one slice of a face including a corner which we do not consider and so there are $\left \lceil{\frac{n}{2}}\right \rceil -1$ distinct clusters or positions where edge cubies can be since the other ones can be reached by a cubie in such a cluster. For example in a $5\times 5 \times 5$ cube, the position next to the corner is one cluster and the one two away from the corner is another. These are the only ones for a $5\times 5 \times 5$ cube since the ones on the other side of the half can be reached from the cubie on the reflection.  Moreover we only need to consider the cubie clusters from one face since the cubies on the other 5 faces also belong to these cubie clusters since they can be reached from a cubie on the face we are considering. Thus since each of the edge clusters can be solved in $O(1)$ moves by \textbf{Lemma 1}, the total number of moves required is $O(n)$ since there are $\left \lceil{\frac{n}{2}}\right \rceil -1$ distinct edge clusters which need to be solved over the whole cube.

We now move on to solving the center cubies. There are $(n-2)^2$ center cubies on one face of a $n\times n \times n$ Rubik's cube since the center cubies form a $n-2$ by $n-2$ square. This leads to $\frac{n^2 - 1}{8}$ unique center cubie clusters when $n$ is odd and $\frac{n^2-2n}{8}$ unique center cubie clusters when $n$ is even. To prove this, we find a recursive formula. When $n$ is odd, we start with a $3\times 3 \times 3$ cube in which case, there is one center cubie and thus one center cluster. Moving on to the $5 \times 5 \times 5$ cube, the center becomes a $3 \times 3$ square on each face and the number of new center cubie clusters is 2 giving a total of 3 center clusters when added to the one center cluster in the $3 \times 3 \times 3$. In general, given the center square of an $n\times n \times n$ cube, the number of new center cubie clusters formed at the $(n+2) \times (n+2) \times (n+2)$ cube is $\left \lceil{\frac{n}{2}}\right \rceil = \frac{n+1}{2}$ since this is the number of edge clusters (see part above) plus the new corner cluster formed on a $n \times n \times n$ cube which becomes a center cubie cluster for the $(n+2) \times (n+2) \times (n+2)$ center portion. Summing over all values for the $n\times n \times n$ cube we see that there are $1+2+ ... +\frac{n-1}{2} = \frac{\frac{n-1}{2} \cdot \frac{n+1}{2}}{2} = \frac{n^2-1}{8}$. When $n$ is even, we start with a $2 \times 2 \times 2$ cube in which case there are 0 center cubies and so 0 clusters. For the $4 \times 4 \times 4$ cube, the center part consists of a $2\times 2 \times 2$ cube which has 1 cubie cluster. Continuing in such a way, going from a $n\times n \times n$ cube to the $(n+2)\times (n+2) \times (n+2)$ cube, we add $\left \lceil{\frac{n}{2}}\right \rceil$ = $\frac{n}{2}$ center cubie clusters by the same reasoning as before. Thus using a similar recurrence, for an $n\times n \times n$ cube, the number of unique center cubie clusters is $1+2+ ... +\frac{n-2}{2} = \frac{\frac{n-2}{2} \cdot \frac{n}{2}}{2} = \frac{n^2-2n}{8}$. Thus for $n$ even or odd, the number of unique center cubie clusters is $O(n^2)$ and since a constant number of moves is required to solve each, the total number of moves required to solve all center cubie clusters is $O(n^2)$. This is the dominant term when summing over edges, corners and centers and so the total number of moves required to solve all clusters in an $n\times n \times n$ Rubik's cube using this naive algorithm is $O(n^2)$.

\section{Optimized Method to Solve the Rubik's Cube using this algorithm}\label{sec:optimizedMethod}
We now present a method to solve the $n\times n \times n$ Rubik's cube in $O(\frac{n^2}{\log{n}})$ moves using this algorithm. The method along with the lemmas have been adapted from \cite{Ref1}.

\textbf{Lemma 2:} Given a solvable $n\times n \times n$ Rubik's configuration, the parity of all cubie clusters can be made even in $O(n)$ moves.

{\textit{Proof.} We assume that the center cluster is already solved since that is usually the first step in solving the Rubik's cube, and therefore we assume that its parity is already even. In addition, any cluster containing at least two indistinguishable cubies can be considered to have even parity or odd parity depending on the chosen label for the indistinguishable cubies. Therefore,we may assume that all such clusters have even parity. This means that all non-edge clusters, including the non-edge cross clusters, can be assumed to have the correct parity no matter how many moves are performed. So we need only fix the parity of the edge clusters.
We start by fixing the parity of the corner cluster and the edge cross cluster (if it exists). Because the cube is solvable, we know that the corner cluster and the edge cross cluster can be solved. Because the corner cluster has $O(1)$ reachable states and the edge cross cluster has $O(1)$ reachable states, we know that we can solve both in $O(1)$ moves. Once those two clusters are solved, we know that their parities must be correct. Therefore, there is a sequence of $O(1)$ moves which can be used to fix the parity of these two clusters. Consider the effect of a face move on the parity of a non-cross edge cluster. For a particular edge cluster, a face move affects the location of eight cubies, due to the fact that a face move also acts like a row or column move for edge cubie groups. The color of each cubie is rotated 90 degrees in the direction of the face's rotation. This means that the permutation applied consists of two permutation cycles each containing four elements. Therefore, if the elements whose colors are changed are $1,2,3,..., 8$ then we can write the applied permutation as $(1 3)(1 5)(1 7)(2 4)(2 6)(2 8)$, or six swaps which means that the parity is still odd (by definition the parity is even if it can be written as the product of an odd number of 2-cycles and odd if it can be written as the product of an even number of 2-cycles). Hence face moves cannot be used to fix the parity of the edge clusters. Now consider the effect of a row or column move on the parity of a non-cross edge cluster. A row or column move affects the colors of four cubies, one for each corner of the rotated slice. The color of each cubie is transferred to the adjacent cubie in the direction of the move rotation. So if the elements whose colors are changed are 1, 2, 3, 4 then the applied permutation is $(1 2 3 4) = (1 2)(1 3)(1 4)$.  Because the permutation can be written as an odd number of swaps, the parity of the cluster has changed. Note, however, that there is exactly one edge cluster whose parity is affected by this movement. Therefore, we can correct the parity of each odd edge cluster by performing a single row or column move that affects the
cluster in question. The total number of moves required is therefore proportional to the number of edge clusters, or $O(n)$. $\Box$

\textbf{Lemma 3:} Suppose we are given an $n \times n \times n$ Rubik's cube configuration and sets $X, Y \subseteq \{0,1,...,\lfloor{\frac{n}{2} - 1} \rfloor \}$ such that $X \cap Y = \emptyset$. If all cubie clusters in $X \times Y$ have the same cluster configuration, then they can all be solved in a sequence of $O(|X| + |Y|)$ moves that only affects cubie clusters in $(X \times Y) \cup (X \times X) \cup (Y \times Y)$.

\textit{Proof.} Let $d$ be the cluster configuration of all clusters in $X \times Y$. By \textbf{Lemma 1}, we know that a constant-length cluster solution exists for $d$. Let $x_1,x_2,..., x_l$ be the elements of $X$ and let $y_1,y_2,...,y_k$ be the elements of $Y$. Then to solve all clusters, the total length of the sequence will be $l + k + 1$. Since the original sequence of moves had length $O(1)$, we know that $l = O(1)$ and so the total length of the solution is $O(|X|+|Y|+1)$. $\Box$

\textbf{Lemma 4:} Suppose we are given an $n \times n \times n$ Rubik's cube configuration, a cluster configuration $c$ and sets $X, Y \subseteq \{0,1,...,\lfloor{\frac{n}{2} - 1} \rfloor \}$ such that $l=|X|$ and $X \cap Y = \emptyset$. Then there exists a sequence of moves of length $O(l\cdot2^l + |Y|)$ such that all cubie clusters $(x,y) \in X \times Y$ in configuration c will be solved, all cubie clusters $(x,y) \in (X \times X) \cup (Y \times Y)$ may or may not be affected and all other cubie clusters will not be affected.

\textit{Proof.} For each row $y \in Y$, let $S_y = \{x \in X |$ cubie cluster $(x,y)$ is in configuration $c$\}. For each set $S \subseteq X$, let $Y_S = \{y \in Y | S_y = S\}$. Because $S \subseteq X$, there are at most $2^l$ different values for $S$. For each S, we will use the results of \textbf{Lemma 3} to construct a sequence of moves to solve each cubie cluster $(x, y) \in S \times Y_S$. This move sequence will have length $O(|S|+|Y_S|)=O(l+|Y_S|)$. Summing over all sets $S \subseteq X$, we get the following number of moves: $O(l\cdot2^l+\sum_{S} |Y_S|) = O(l\cdot 2^l + |Y|)$. $\Box$

\textbf{Lemma 5:} Suppose we are given an $n \times n \times n$ Rubik's Cube configuration, a cluster configuration $c$, and sets $X,Y \subseteq \{0,1,...,\lfloor{\frac{n}{2} - 1} \rfloor \} $  such that $X \cap Y = \emptyset$. Then there exists a sequence of moves of length $O(\frac{|X| \cdot |Y|}{\log |Y|})$ such that all cubie clusters $(x,y) \in X \times Y$ in configuration $c$ wil be solved, all cubie clusters $(x,y) \in (X \times X) \cup (Y\times Y)$ mar or may not be affected and all other cubie clusters will not be affected.

\textit{Proof.}  Let $l = \frac{1}{2} \cdot \log_2 |Y|$, such that $2^l = \sqrt{|Y|}$. Let $k = \left \lceil{\frac{|X|}{l}} \right \rceil$. Partition the set $X$ into a series of sets $X_1, ..., X_k$ each of which has size $\leq l$. For each $X_i$, we solve the cubie clusters in $X_i \times Y$ using the sequence of moves that is guaranteed to exist by \textbf{Lemma 4}. The number of moves required to solve a single $X_i$ is:
\begin{align*}
     O(l\cdot 2^l + |Y|) = O((\frac{1}{2}\cdot \log_2 |Y|) \cdot \sqrt{|Y|} + |Y|) = O(|Y|)\\
\end{align*}
\normalsize{Therefore if we wish to perform this for $k$ sets, the total number of moves becomes}
\begin{align*}
    O(k \cdot |Y|) = O (\frac{|X| \cdot |Y|}{\frac{1}{2} \cdot \log_2 |Y|}) = O(\frac{|X| \cdot |Y|}{\log |Y|}).\\
\end{align*}
$\Box$

\textbf{Theorem:} Given an $n \times n \times n$ Rubik's cube configuration, all cubie clusters can be solved in $O(\frac{n^2}{\log n})$ moves.

\textit{Proof.} In order to solve the Rubik's cube, we begin by fixing the parity which we know by \textbf{Lemma 2} can be done in $O(n)$ moves. Then we solve each edge cluster, each of which can by \textbf{Lemma 1} be solved in $O(1)$ moves, so this step takes $O(n)$ moves. Once the edges have been solved, we solve the non-edge clusters. Let $k = \sqrt{\frac{n}{2}}$. Partition $\{0,1,...,\lfloor{\frac{n}{2}} \rfloor - 1\}$ into sets $G_1, G_2,..., G_k$ each with size $ \leq k$. For each pair $i,j$ such that $i \neq j$ and each cluster configuration $c$, we use the sequence of moves guaranteed by \textbf{Lemma 5} to solve all $(x,y) \in G_i \times G_j$ with the configuration $c$. This ensures that all cubie clusters $(x,y) \in G_i \times G_j$ will be solved. For each $i$, we must also solve all cubie clusters $(x,y) \in G_i \times G_i$. There are $k^2\cdot \frac{k-1}{2}$ such cubie clusters, and we solve each one individually. For a single pair $i \neq j$ and a single configuration $c$, the number of moves required will be $O(\frac{k^2}{\log k}) = O(\frac{n}{\log n})$. There are a constant number of possible configurations so solving a single pair $i,h$ for all configurations will also require $O(\frac{n}{\log n})$ moves. There are $k^2 - k$ such pairs which thus gives a total of $O(\frac{n^2}{\log n})$ moves as desired. $\Box$

\subsection{Comparison of Algorithms}\label{ssec:comparison}

Using the method in \cite{Ref1}, we have shown that our algorithm can be used to solve the $n \times n \times n$ cube in $O(\frac{n^2}{\log n})$. We can compare the performance of our algorithm to the algorithm presented in \cite{Ref1}, as shown in Table \ref{table:table1}.

\begin{table}[h!]
\caption{Comparison of algorithms}
\label{table:table1}
  \begin{center}
    \begin{tabular}{l|c|r} 
      \textbf{Algorithm} & \textbf{Complexity} & \textbf{Number of Moves Per Permutation}\\
      \hline
      Algorithm proposed in this paper & $O(\frac{n^2}{\log n})$ & $8$ or $10$\\
      Algorithm presented in \cite{Ref1} & $O(\frac{n^2}{\log n})$ & $10$\\
    \end{tabular}
  \end{center}
\end{table}

Our algorithm improves the results of \cite{Ref1} in terms of the total number of moves required since our algorithm uses 8 or 10 moves per permutation while the algorithm in \cite{Ref1} uses 10 moves for every permutation. The algorithmic complexity of both algorithms using the method is the same at $O(\frac{n^2}{\log n})$.
\section{Conclusion}\label{sec:conclusion}

In this paper, we presented three possible variations of an algorithm, each of minimum length 8 and maximum length 10, for permuting corner, edge and center cubies in any $n \times n \times n$ Rubik's cube and used this to show that the whole $n \times n \times n$ cube can be solved in $O(n^2)$ moves. We then applied the method from \cite{Ref1} to reduce the complexity of the algorithm to $O(\frac{n^2}{\log n})$. Using our algorithm with the method from \cite{Ref1} also reduced the total number of moves required to solve the cube compared to the algorithm used in \cite{Ref1}.

\ifCLASSOPTIONcaptionsoff
  \newpage
\fi


\begin{thebibliography}{1}

\bibitem{Ref1}
Erik Demaine, Martin Demaine, Sarah Eisenstat, Anna Lubiw and Andrew Winslow, \emph{Algorithms for Solving Rubik's Cubes}, https://arxiv.org/pdf/1106.5736.pdf


\end{thebibliography}
\end{document}